\documentclass[a4paper,11pt]{amsart}
\usepackage{amsfonts,amssymb,amsthm,cite,amsmath,amstext}
\usepackage{color}
\usepackage[colorlinks,linkcolor=black,citecolor=black]{hyperref}
\usepackage{comment}
\usepackage{framed}
\usepackage{calligra} 
\usepackage{tikz-cd}

\definecolor{shadecolor}{gray}{0.875}

\newtheorem{thrm}{Theorem}[section]
\newtheorem{thrmx}{Theorem}

\newtheorem{corx}{Corollary}

\newtheorem{lem}[thrm]{Lemma}
\newtheorem{cor}[thrm]{Corollary}
\newtheorem{prop}[thrm]{Proposition}
\newtheorem{conj}[thrm]{Conjecture}

\theoremstyle{definition}
\newtheorem{defn}[thrm]{Definition}

\newtheorem{rmk}[thrm]{Remark}

\newenvironment{claim}
            {\par \bigskip \noindent \textbf{Claim:}}
            {$\Box$ \par \bigskip \noindent}

\DeclareMathOperator{\relint}{relint}

\DeclareMathOperator{\vol}{vol}
\DeclareMathOperator{\Mov}{Mov}

\DeclareMathOperator{\supp}{supp}

\DeclareMathOperator{\Image}{Image}

\DeclareMathOperator{\nd}{nd}

\DeclareMathOperator{\Span}{span}

\DeclareMathOperator{\codim}{codim}

\DeclareMathOperator{\prim}{Prime}

\DeclareMathOperator{\Nef}{Nef}
\DeclareMathOperator{\Psef}{Psef}

\DeclareMathOperator{\Null}{Null}
\DeclareMathOperator{\eff}{eff}
\DeclareMathOperator{\zar}{zar}

\title
[Numerical characterization of the hard Lefschetz classes of dimension two]{Numerical characterization of the hard Lefschetz classes of dimension two, II:
\\
\footnotesize{supercritical collections of free divisor classes}}
\author{Jiajun Hu and Jian Xiao}
\date{}

\begin{document}


\begin{abstract}
For $(n-2)$ free divisor classes on a smooth projective variety of dimension $n$, the product of these free divisor classes induces a Lefschetz type operator acting on the N\'{e}ron-Severi space or the cohomology group of $(1,1)$ classes.
We give a characterization of this kernel space, when the collection of these free divisor classes is supercritical. This resolves Shenfeld-van Handel's open problem in this setting. As consequences, we provide an algebro-geometric proof of the characterization of the extremals of the Alexandrov-Fenchel inequality for a supercritical collection of rational convex polytopes; we also give a characterization of the extremals of the Khovanskii-Teissier inequality given by the intersection numbers of two arbitrary free divisor classes.

\end{abstract}

\maketitle

\tableofcontents

\section{Introduction}

In this article, we work on projective varieties over the field of complex numbers.

\subsection{Motivation}
Let $K_1,...,K_m$ be convex bodies in $\mathbb{R}^n$ and let $t_1,...,t_m$ be nonnegative real numbers, a basic result due to Minkowski tells that the volume function $\vol(\sum_{i=1}^m t_i K_i)$ is a homogeneous polynomial of $t_1,...,t_m$:
\begin{equation*}
\vol(\sum_{i=1}^m t_i K_i) =\sum_{i_1+...+i_m =n} V(K_{i_1},...,K_{i_m}) t_{i_1}...t_{i_m}.
\end{equation*}
The coefficients $V(K_{i_1},...,K_{i_m})$ are called mixed volumes. For the basics on mixed volumes, we refer the reader to Schneider's book \cite{schneiderBrunnMbook}. By taking certain tuple of convex bodies $K_1,...,K_m$, one can recover many geometric quantities of a convex body, such as volume, projection volume, surface area, mean width, or more generally, mean size in intermediate dimension -- quermassintegral.  
A cornerstone in convex geometry is the Alexandrov-Fenchel (AF for short) inequality: for convex bodies $K, L$ and a collection of $(n-2)$ convex bodies $\mathcal{C}=(C_1,...,C_{n-2})$ in $ \mathbb{R}^n$, we have
\begin{equation*}
  V(K, L, \mathcal{C})^2 \geq V(K, K, \mathcal{C}) V(L, L, \mathcal{C}),
\end{equation*}
where $V(K_1, K_2,\mathcal{C})$ is the mixed volume of the convex bodies $K_1, K_2, C_1,...,C_{n-2}$. This result has profound impacts in convex geometry. Many geometric inequalities, such as the isoperimetric inequalities and the Brunn-Minkowski inequalities, are consequences of the AF inequality. 
A foundational question is the characterization of its extremals, that is, 
\begin{center}
  (\textbf{AF}) \emph{When does the equality hold in the AF inequality?}
\end{center}
More precisely, if we fix the collection $\mathcal{C}$, what is the relation between $K, L$ and the collection $\mathcal{C}$ when equality holds?
This is a long-standing open problem, which dates back to Alexandrov's original paper \cite{alexandroff1938theorie}.
In the recent breakthrough work \cite{handelAFextremACTA}, Shenfeld and van Handel completely solved the (\textbf{AF}) problem when $K, L$ and the collection $\mathcal{C}$ are convex polytopes. For the history, some previous development, results on convex bodies beyond polytopes and the difficulty in resolving this fundamental problem, we refer the reader to the exposition in \emph{loc. cit.} and the references therein.

The AF inequality goes far beyond convexity, it has various counterparts and applications in other areas of mathematics (see e.g. \cite{stanleyAF, stanleyposet, stanleyLogConc, huh2022combinatorics} and the references therein). In particular, regarding its interactions with Hodge theory, in algebraic geometry, its counterpart is the Hodge index theorem or Khovanskii-Teissier inequalities (see e.g. \cite{lazarsfeldPosI}); in combinatorics, we have its analogs in combinatorial Hodge theory \cite{huh2022combinatorics}.
Due to the deep connections of the AF inequality with other areas of mathematics, and the expectation that the algebraic formulation behind the AF inequality may carry over to other mathematical problems, Shenfeld-van Handel \cite{handelAFextremACTA} suggested an open proposal to establish the algebraic analogues of the characterization of the extremals. The central topic of this article is around Shenfeld-van Handel's open problem.

This article is also motivated by another closely related (and more general) question in algebraic geometry. 
Let $X$ be a smooth projective variety of dimension $n$. We denote by $H^{p,p}(X)$ the cohomology group of $(p,p)$ classes with real coefficients. 
The vector space of numerically equivalent cycle classes of codimension $p$ (equivalently, of dimension $n-p$) with real coefficients is denoted by $N^p (X)$ (equivalently, $N_{n-p} (X)$). The space $N^1 (X)$ is known as the N\'eron-Severi space with real coefficients. A divisor class $L$ is called nef iff $L\cdot [C]\geq 0$ for any irreducible curve $C$ on $X$. 
Fix a collection of nef classes $\mathcal{L}=\{L_1, ...,L_m\}$ on $X$, and denote the complete intersection class by $$\mathbb{L}=L_1 \cdot ...\cdot L_m.$$ The class $\mathbb{L}$ can be viewed as a linear operator acting on some cohomology group $H^{p,p}(X)$ or some cycle classes group $N^p (X)$. A natural question is, 
\begin{quote}
\begin{center}
 (\textbf{Ker}) \emph{What is the kernel space of $\mathbb{L}$?}
\end{center}
\end{quote}
The (\textbf{Ker}) problem has rich connection with the geometry of algebraic cycles and algebraic maps. A complete resolution to this question will be very challenging, as even in some special cases it could imply a strong form of the generalized Hodge conjecture for smooth projective varieties \cite{voisinPushforw}.  Nevertheless, the question is accessible when $m=n-2$, in which case
the class $\mathbb{L}$ induces the following Lefschetz-type operator
\begin{equation*}
  \mathbb{L}: H^{1,1} (X) \rightarrow H^{n-1,n-1} (X) .
\end{equation*}
It is sure that one may also consider the operator acting on the N\'eron-Severi space $N^1 (X)$. One intriguing aspect on the (\textbf{Ker}) problem is that, 
the algebro-geometric counterpart of the characterization of the extremals for (\textbf{AF}) is exactly the (\textbf{Ker}) problem when we restrict on the case when $\mathcal{L}$ is a collection of $(n-2)$ classes:
 \begin{center}
(\textbf{Ker$_2$})  \emph{Characterize the kernel space of $\mathbb{L}: H^{1,1} (X) \rightarrow H^{n-1,n-1} (X)$.}
\end{center}
We shall give a complete solution to the (\textbf{Ker$_2$}) problem when $\mathcal{L}$ is a supercritical collection of free classes. This resolves Shenfeld-van Handel's open problem in this setting.

\subsection{The conjecture}
We first introduce some notations. 
For a positive integer $m$, $[m]$ denotes the finite set $\{1,2,...,m\}$. Given finite vectors $v_1,...,v_m$ in a vector space $V$ and a non-empty subset $I\subset [m]$, we denote
$$v_I = \sum_{i\in I} v_i.$$
Given a real vector space $V$ and a subset $S \subset V$, the linear subspace generated by $S$ is denoted by $\Span_\mathbb{R} (S)$. We denote the set of prime divisors on $X$ by $\prim(X)$. The numerical dimension of a nef class $L$ is defined by $$\nd(L)=\max\{k \geq 0 | L^k \neq 0\}.$$
For a collection of nef classes $\mathcal{L}=\{L_1, ...,L_{n-2}\}$ on $X$, the supercritical condition says that
\begin{equation*}
  \nd(L_I) \geq |I|+2, \forall I\subset[n-2].
\end{equation*}

Now we can state Shenfeld-van Handel's conjecture\footnote{While the paper \cite[Section 16]{handelAFextremACTA} did not explicitly state the conjecture, it is not hard to formulate the conjecture on a toric variety by their characterization on the extremals for convex polytopes and then extend the conjecture to an arbitrary projective variety via positivity theory in algebraic geometry.} for a supercritical collection:

\begin{conj}\label{main conj}
Let $X$ be a smooth projective variety of dimension $n$. Let $\mathcal{L}=(L_1,...,L_{n-2})$ be a supercritical collection of nef classes on $X$. Define the vector spaces $V_{\mathcal{L},\eff}$ as follows:
\begin{align*}
V_{\mathcal{L}, \eff} = \Span_{\mathbb{R}} \{[D]: D\in \prim(X),\ \mathbb{L} \cdot [D]=0\}.
\end{align*}
Then $$\ker \mathbb{L} = V_{\mathcal{L},\eff}.$$
In particular, the class $\mathbb{L}$ is a hard Lefschetz class if and only if $\mathbb{L} \cdot [D] \neq 0$ for any prime divisor $D$ on $X$.
\end{conj}

When the collection is not supercritical, it is not hard to find examples where the space $V_{\mathcal{L},\eff}$ is not sufficient to characterize $\ker \mathbb{L}$. One then needs another subtler piece, called the ``degenerate'' part. For the precise conjectural picture and related discussions in this more general setting, see \cite[Conjecture 1.9 and Section 7]{huxiaoLefcharArxiv}.

\subsection{Previous works}

We list several previous works that are closely related to the (\textbf{Ker$_2$}) problem and our contexts.
\begin{enumerate}
  \item The first important result is due to Panov \cite{Panov1987ONSP}, who found the dimensionality mechanism in Alexandrov's mixed discriminant inequality. Panov's result solves the linear-algebra version of (\textbf{Ker$_2$}). It can be reformulated in a more algebro-geometric scheme on a compact complex torus 
      (see \cite[Section 7]{huxiaoLefcharArxiv}). 
  \item A complete solution to (\textbf{Ker$_2$}) in the setting of arbitrary convex polytopes is due to Shenfeld-van Handel \cite{handelAFextremACTA}. There is no doubt that this work paves the way to the conjectural picture in the algebraic setting, and leads to an expected answer to (\textbf{Ker$_2$}).
  \item In \cite{huxiaoLefcharArxiv}, we solved the (\textbf{Ker$_2$}) problem when the nef collection $\mathcal{L}=(L_1,...,L_{n-2})$ satisfies the positivity condition: $$\nd(L_i)\geq i+2 \text{ for any } i\in [n-2].$$ 
      In particular, this includes the case of big nef collections.
      Though every supercritical collection has a rearrangement satisfying this stronger requirement, it is not quite clear yet how to relate the solution to a collection under rearrangement with a general supercritical condition \cite[Proposition 6.8]{huxiaoLefcharArxiv}.
  \item Beyond (\textbf{Ker$_2$}), in our previous works \cite{huxiaohardlef2022Arxiv} and \cite{hushangxiao2023hardlefArxiv}, we initiate the program on the numerical characterization problem for hard Lefschetz classes (of any degrees) arising from weak positivity on a projective variety (or in the analytic setting, on a compact K\"ahler manifold), which is closely related to the (\textbf{Ker$_d$}) problem. See also \cite{decataldoLefsemismall} for some motivation from its relation with Hodge theory of algebraic maps. In particular, we obtained some sufficient positivity criterion on $\ker \mathbb{L} =0$ when $\mathcal{L}$ is a collection of $(n-2d)$ classes and $\mathbb{L}$ is the Lefschetz type operator acting on $H^{d,d}(X)$:
      \begin{equation*}
  \mathbb{L}: H^{d,d} (X)\rightarrow H^{n-d,n-d} (X).
\end{equation*} 
The attempt to understand the sufficient criterion in more depth leads us to realize that the numerical characterization of hard Lefschetz classes of dimension two is also closely related to Shenfeld-van Handel's open proposal.
\end{enumerate}

\subsection{Main results}
As mentioned above, in \cite{huxiaoLefcharArxiv}, we proved Conjecture \ref{main conj} when the collection $\mathcal{L}=(L_1,...,L_{n-2})$ satisfies that $\nd(L_i)\geq i+2$ for any $i\in [n-2]$. Though this result only needs nefness, this positivity requirement is not quite satisfactory, as it is stronger than supercriticality.
In this paper, we confirm Conjecture \ref{main conj} in the free case, with a self-contained proof.

For a divisor class $L$ with $\mathbb{Z}$-coefficients, we will not distinguish the divisor class and the associated line bundle, which we still denote by the same symbol $L$. The divisor class/line bundle $L$ is called free if the linear system $|L|$ is basepoint-free, or equivalently, the Kodaira map given by $H^0 (X, L)$
\begin{align*}
  \varphi_L : X \rightarrow \mathbb{P}(H^0 (X, L)),
  \ x\mapsto [s_0 (x):\cdots:s_N (x)] 
\end{align*}
is a morphism, where $\{s_k\}$ is a basis of $H^0 (X, L)$. When $L$ is free, it is clear that the numerical dimension $\nd(L)=\dim Y_L$ where $Y_L$ is the image of $\varphi_L$.

\begin{thrmx}\label{mainthrm}
Conjecture \ref{main conj} holds when $\mathcal{L}=(L_1,...,L_{n-2})$ is a supercritical collection of free divisor classes.
\end{thrmx}

For the collection in the setting of Theorem \ref{mainthrm} and the prime divisor $D$ appearing in $ V_{\mathcal{L},\eff}$, one can prove that the number of such $D$ is at most the Picard number of $X$, and these divisors are contained in the augmented base locus of a big class and have very strong rigidity, and in particular, span extremal rays of the pseudo-effective cone (see Remark \ref{rigid}).

\begin{rmk}
Recall that a divisor class $L$ is called semiample if $kL$ is free from some positive integer $k$, therefore Theorem \ref{mainthrm} also holds when $\mathcal{L}=(L_1,...,L_{n-2})$ is a supercritical collection of semiample divisor classes.
\end{rmk}

Comparing to our previous result \cite[Theorem A]{huxiaoLefcharArxiv}, though with an extra requirement on freeness, the advantage of Theorem \ref{mainthrm} is clear -- it works for an arbitrary supercritical collection. This freeness setting is already very useful. As a main application, we obtain an algebro-geometric proof of the characterization of the extremals of the AF inequality for a supercritical collection of rational convex polytopes, originally due to \cite[Corollary 2.16 or Theorem 8.1]{handelAFextremACTA}:

\begin{corx}\label{corhandel}
Let $\mathcal{P}=(P_1,...,P_{n-2}) \subset \mathbb{R}^n$ be a supercritical collection of rational polytopes. Assume that $K,L$ are rational polytopes such that $V(K, L,\mathcal{P} )>0$. Then 
\begin{equation*}
  V(K, L, \mathcal{P})^2 = V(K, K, \mathcal{P})V(L, L, \mathcal{P})
\end{equation*}
if and only if there exist $a>0$ and $v\in \mathbb{R}^n$ such that $K$ and $aL+v$ have the same
supporting hyperplanes in all $(B, \mathcal{P})$-extreme normal directions. 
\end{corx}

Further terms on convexity in Corollary \ref{corhandel} can be found in Section \ref{sec appl}. When $V(K, L,\mathcal{P} )=0$, by the AF equality we automatically have equality. The characterization of $V(K, L,\mathcal{P} )=0$ is provided by the dimensionality mechanism. 
Corollary \ref{corhandel} follows from Theorem \ref{mainthrm}, the construction via toric geometry and the fact that on a smooth projective toric variety, every nef line bundle is semiample.

Note that many log-concave phenomenon, frequently appearing in combinatorics, can be reformulated by the Khovanskii-Teissier inequality on the intersection numbers of two free divisor classes. Our another application is the characterization of the extremals of the Khovanskii-Teissier inequality given by the intersection numbers of two free divisor classes. 
 
 \begin{corx}\label{corlogc}
Let $X$ be a smooth projective variety of dimension $n$, and let $A, B$ be free divisor classes on $X$. Fix $1\leq k\leq n-1$. Assume that $(A^k \cdot B^{n-k})>0$,  then 
\begin{equation*}
  (A^k \cdot B^{n-k})^2 =  (A^{k-1} \cdot B^{n-k+1}) (A^{k+1} \cdot B^{n-k-1})
\end{equation*}
if and only if there exist some $c>0$ and $c_i\in \mathbb{R}$ such that
\begin{equation*}
  A-cB= \sum_{i=1}^N c_i [D_i]
\end{equation*}
where the $D_i$ are primes divisors such that $A^{k-1} \cdot B^{n-k-1}\cdot [D_i] =0$.
\end{corx}

When $(A^k \cdot B^{n-k})=0$, the equality automatically holds, in which case by using \cite{huxiaohardlef2022Arxiv}, we have $(A^k \cdot B^{n-k})=0$ iff either $\nd(A)<k$, or $\nd(B)<n-k$, or $\nd(A+B)<n$.

This article is organized as follows. In Section \ref{sec pre} we recall some basic notions and present several preliminary results. Section \ref{sec main result} is devoted to the proof of the main theorem. We begin with the baby case on 3-folds, and as an instructive prototype we extend the baby case to higher dimensions by showing how the induction on dimension works, and finally we apply these results to prove the general case. In Section \ref{sec appl}, we provide an algebro-geometric proof of the characterization of the extremals of the Alexandrov-Fenchel inequality for a supercritical collection of rational convex polytopes; we also give a characterization of the extremals of the Khovanskii-Teissier inequality given by the intersection numbers of two arbitrary free divisor classes.

\section{Preliminaries}\label{sec pre}

In this section, we introduce some basic notions and preliminary results which will be applied in the sequel.

While we work on a projective variety $X$ defined over the field of complex numbers and study the kernel spaces in $H^{1,1}(X)$. All statements remain valid over any algebraically closed field of characteristic zero upon replacing the cohomology group $H^{1,1}(X)$ by the N\'eron-Severi space $N^1(X)$. We use the following convention.
For a subvariety  $Z\subset X$, we denote by $$[Z]_{X,hom}\in H_{2\dim Z}(X,\mathbb{R})$$ the fundamental class of $Z$ in the homology group. When $X$ is smooth, we write
$$[Z]_X\in H^{2\dim X-2\dim Z}(X,\mathbb{R})$$ 
for the Poincar\'e dual of $[Z]_{X,hom}$, omitting the ambient subscript $X$ if the context is clear. For subvarieties $Z\subset Y\subset X$, we write $j_{Z\subset Y}$ for the closed embedding $$j_{Z\subset Y}:Z\hookrightarrow Y.$$ 
    In the special case $Y=X$, we abbreviate $j_{Z\subset X}$ simply as $j_Z$.

\subsection{Positivity}
We recall several positivity notions on a projective variety. For more notions on positivity in algebraic geometry and analytic geometry, see \cite{dem_analyticAG, lazarsfeldPosI}.

\begin{defn}
Let $X$ be a smooth projective variety of dimension $n$.
\begin{itemize}
\item $N^1(X)$ is the real vector space of numerically equivalent classes of divisors.
\item $N_1(X)$ is the real vector space of numerically equivalent classes of curves.
\item $\prim(X)$ is the set of prime divisors on $X$.
\item $\Psef^1(X)$ is the cone of pseudo-effective divisor classes, that is, any element $L\in \Psef^1(X)$ can be written as $$L=\lim_k [F_k],$$ where $F_k$ is a non-negative combination of elements in $\prim(X)$. Then $\Psef^1(X)$ is a full-dimensional closed convex cone in $N^1(X)$. The interior of $\Psef^1(X)$ is called the big cone, and an interior point is called a big class.
\item $\Nef^1(X) $ is the cone of nef divisor classes, whose interior is the ample cone. 
\end{itemize}

In the above algebraic setting, all the cones are either in $N^1(X)$ or $N_1(X)$.

\end{defn}

More generally, the algebraic positivity notions have their analytic counterparts when one studies the K\"ahler geometry of a smooth projective variety.

\begin{defn}
Let $X$ be a compact K\"ahler manifold of dimension $n$.
\begin{itemize}
\item $H^{1,1}(X)$ is the real vector space of $(1,1)$ classes.
\item $H^{n-1,n-1}(X)$ is the real vector space of $(n-1,n-1)$ classes.
\item $\Psef^1(X)$ is the cone of pseudo-effective $(1,1)$ classes, that is, $L \in \Psef^1(X)$ if and only if it can be represented by a positive $(1,1)$ current. An interior point of this cone is called a big $(1,1)$ class.
\item $\Nef^1(X)$ is the cone of nef $(1,1)$ classes. It is the closure of the K\"ahler cone, generated by K\"ahler classes.
\end{itemize}

\end{defn}

On a smooth projective variety, $\alpha\in \Nef^1(X)$  iff for any irreducible curve $C$, $\alpha\cdot [C]\geq 0$. For a nef class $\alpha\in \Nef^1(X)$, its numerical dimension is defined by $$\nd(\alpha)=\max\{k \geq 0 | \alpha^k \neq 0\}.$$
Then a nef class $\alpha\in \Nef^1(X)$ is big iff the numerical dimension $\nd(\alpha)=n$. 
In the analytic setting, we shall need to consider the movable cone $\Mov^1(X)$ of $(1,1)$ classes:  $L\in \Mov^1(X)$ if and only if $$L = \lim_{k\rightarrow \infty} (\pi_k)_* \omega_k,$$ where $\pi_k: X_k \rightarrow X$ is a K\"ahler modification and $\omega_k$ is a K\"ahler class on $X_k$. It is clear that 
\begin{equation*}
  \Nef^1(X) \subset \Mov^1(X) \subset \Psef^1(X), 
\end{equation*}
and in general the inclusions are strict. A remarkable result relating the cones $\Psef^1(X), \Mov^1(X)$ is the divisorial Zariski decomposition due to Boucksom \cite{Bou04}.

\begin{lem}
\label{zariski}
Every $L\in \Psef^1(X)$ can be decomposed as follows:
\begin{equation*}
  L=P(L)+[N(L)],
\end{equation*}
where $P(L) \in \Mov^1(X)$ is called the positive part, and $N(L)$ is an effective divisor which is called the negative part. 
 Furthermore, $N(L)$ is supported by at most $\rho(X)$ prime divisors, where $\rho(X)$ is the Picard number of $X$.
\end{lem}

\subsection{Proportionality via Lorentzian signature}

The following result is well known (see e.g. \cite[Proposition 3.4]{huxiaoLefcharArxiv}).

\begin{prop}\label{proportional}
Let $(V,Q)$ be a real vector space $V$ of finite dimension, together with a quadratic form $Q$ with exactly one positive eigenvalue. Assume that $x, y\in V$ satisfies that $Q(x)\geq 0, Q(y)\geq 0$ and
$$Q(x, y)^2 = Q(x)Q(y),$$
then the linear forms $Q(x, -)$ and $Q(y, -)$ are proportional.
\end{prop}

As a consequence, we obtain:

\begin{lem}\label{equalAF}
    Let $\mathcal{L}=\{L_1,...,L_{n-2}\}$ be a collection of nef classes, and let $\alpha,\beta \in H^{1,1}(X,\mathbb{R})$ satisfy
    \begin{equation*}
        (\mathbb{L} \cdot \alpha\cdot\beta)^2 =(\mathbb{L} \cdot \alpha^2)(\mathbb{L} \cdot \beta^2) \text{ and } \mathbb{L} \cdot\alpha^2 \geq 0, \mathbb{L} \cdot\beta^2 \geq 0
    \end{equation*}
    then $\mathbb{L}\cdot \alpha$ and $\mathbb{L}\cdot \beta$ are proportional.
    Moreover, if 
    \begin{equation*}
       \mathbb{L} \cdot \alpha\cdot\beta =\mathbb{L} \cdot \alpha^2=0 \text{ and } \mathbb{L} \cdot\beta^2 >0,
    \end{equation*}
     then $\mathbb{L} \cdot\alpha=0$.
\end{lem}

\begin{proof}
This follows immediately from Proposition \ref{proportional} by setting $Q(\alpha, \beta)=\mathbb{L} \cdot\alpha\cdot\beta$ and using Hodge index theorem (i.e., Hodge-Riemann relation in degree one).
\end{proof}

\subsection{Hall-Rado for nef classes}

\begin{defn}
    Given a tuple of nef classes $\mathcal{L}=(L_1,...,L_m)$, we put $\mathbb{L}=L_1\cdot...\cdot L_m$ and define the numerical dimension of the collection $\mathcal{L}$ as
    \begin{equation*}
        \nd(\mathcal{L})=\min_{\emptyset \neq I\subset [m]}\{\nd(L_I)+m-|I|\}.
    \end{equation*}
    It is called supercritical if $\nd(\mathcal{L})\geq |\mathcal{L}|+2$.
\end{defn}

\begin{lem}\label{hallrado}
    Let $\mathcal{L}=(L_1,...,L_m)$ be a collection of nef classes with $m\leq n$. Then 
    $$\nd(\mathcal{L})\geq m \iff \mathbb{L}\neq 0.$$
\end{lem}

\begin{proof}
Note that by definition $\nd(\mathcal{L})\geq m$ iff $\nd(L_I)\geq |I|$ for any $I\subset [m]$. 

The result follows from the Hall-Rado relation for nef classes proved in \cite{huxiaohardlef2022Arxiv} (see also \cite[Theorem 3.7]{huxiaoLefcharArxiv} for the extension to Lorentzian polynomials).
\end{proof}

\subsection{Null locus of a collection of classes}

\begin{defn}\label{locusdef}
Let $X$ be a smooth projective variety of dimension $n$ and let $\mathcal{L}=(L_1,...,L_{n-2})$ be a collection of nef classes. The null locus of codimension $l$ of $\mathcal{L}$ is defined by
\begin{equation*}
  \Null^l (\mathcal{L})= \overline{\bigcup_{\mathbb{L}\cdot [V]=0, \codim V=l} V}^{\zar}
  \end{equation*}
where $V$ ranges over irreducible subvarieties of codimension $l$ with $\mathbb{L}\cdot[V]=0$. 
\end{defn}

A key property is that the supercritical condition ensures the properness of $\Null^2 (\mathcal{L})$. This has been noted in \cite{huxiaoLefcharArxiv}. For completeness and reader's convenience, we include it below.

\begin{lem}\label{locus}
Let $X$ be a smooth projective variety of dimension $n$ and let $\mathcal{L}=(L_1,\cdots,L_{n-2})$ be a supercritical collection of free divisor classes.
Then $\Null^2 (\mathcal{L})$ is a proper Zariski closed subset.
\end{lem}

\begin{proof}
Since the collection $\mathcal{L}=(L_1,\cdots,L_{n-2})$ is free, for any $I\subset [n-2]$, the line bundle $L_I$ is free.
Let $$\phi_I:X \rightarrow \mathbb{P}(H^0(X,L_I))$$ be the Kodaira map induced by $L_I$ and let $Y_I=\phi_I(X)$ be the image of $\phi_I$.
Let $V \subset X$ be a codimension-two subvariety with $\mathbb{L}\cdot [V]=0$. By Lemma \ref{hallrado}, there exits $I\subset [n-2]$ such that
\begin{equation*}
    L_{I}^{|I|}\cdot [V]=0.
\end{equation*}
This is equivalent to that $\dim \phi_I(V) <|I|$,
which implies
  \begin{equation*}
    V \subset \phi_I^{-1}(Y_I^{\geq n-|I|-1})
  \end{equation*}
where
  \begin{equation*}
    Y_I^{\geq n-|I|-1}:=\{y\in Y_I: \dim \phi_I^{-1}(y)\geq n-|I|-1\}.
  \end{equation*}
Since $V$ is arbitrary, we obtain
\begin{equation*}
    \Null^2(\mathcal{L})\subset \bigcup_{I\subset [n-2]}\phi_I^{-1}(Y_I^{\geq n-|I|-1}).
\end{equation*}
It remains to show the properness of each $\phi_I^{-1}(Y_I^{\geq n-|I|-1})$. To this end, we use the supercriticality of the collection $\mathcal{L}$ which implies that
\begin{equation*}
  \dim Y_I = \nd(L_I)\geq  |I|+2.
\end{equation*}
Hence, a general fiber of $\phi_I$ has dimension at most $n-|I|-2$ and as a consequence, $$Y_I^{\geq n-|I|-1}\subset Y_I$$ is a proper closed subset. This finishes the proof.
\end{proof}

\begin{rmk}
The properness in Lemma \ref{locus} still holds if the supercritical collection $\mathcal{L}$ is only nef, this general case is proved in \cite[Lemma 6.4]{huxiaoLefcharArxiv}, where we applied a deep result on the properness of the null locus of a big nef class (see e.g. \cite{tosattiNullLocus}). 
\end{rmk}

\begin{rmk}\label{rigid}
Assume that the nef collection $\mathcal{L}=(L_1,...,L_{n-2})$ satisfies the condition that for any subset $I\subset [n-2]$, the numerical dimension satisfies $\nd(L_I)\geq |I|+1$. Then the number of prime divisors $D$ such that $\mathbb{L}\cdot [D]=0$ is at most the Picard number of $X$. Moreover, each such $D$ has strong rigidity: it appears in the negative part of some big class (hence, it lies in the augmented base locus of the big class), and in particular, spans an extremal ray of the pseudo-effective cone $\Psef^1(X)$.
This result follows from \cite[Proposition 5.6]{huxiaoLefcharArxiv}.
\end{rmk}

\subsection{Homology class of a general fiber}
The following recent result due to \cite[Lemma 2.6]{huh2025arxiv-realizationshomologyclassesprojection} is helpful, where they applied it to a different context.

\begin{lem}\label{fiberclass}
Let $f: X\rightarrow Y$ be a surjective morphism between irreducible projective varieties over an algebraically closed field, then the irreducible components of a general fiber of $f$ are algebraically equivalent to each other in $X$. 
\end{lem}

As a consequence, we obtain:
\begin{lem}\label{algequiv}
Let $X$ be a projective variety of dimension $n$. Let $D$ be a prime divisor and let $L$ be basepoint-free divisor which is general in its linear series. Then for any irreducible component $W$ of $D\cap L$, there exists $c_W\in \mathbb{R}$ such that
    $$[W]_{D,hom} = c_W [L]_{|D}\cap [D]_{D,hom} \text{ in } H_{2n-4}(D,\mathbb{R}).$$
\end{lem}
\begin{proof}
    Let $$\phi:X\rightarrow Y\subset \mathbb{P}(H^0(X,L))$$ be the Kodaira map defined by the line bundle $L$. Then $L=\phi^*A$ where $A=\mathcal{O}(1)_{|Y}$ is the hyperplane line bundle on $Y$.
     
     If $\dim \phi(D)=0$, then for a general divisor $H\in |L|$, $H\cap D=\emptyset$, so there is nothing to prove. If $\dim \phi(D)\geq 2$, then by Bertini's theorem, the intersection  $H\cap D$ is irreducible 
    for general $H\in |L|$ and the result follows as well. Now suppose $\dim \phi(D)=1$, so that $C=\phi (D)\subset Y$ is a curve. Consider the restriction map
    \begin{equation*}
        f=\phi_{|D}:D\rightarrow C.
    \end{equation*}
    Applying Lemma \ref{fiberclass} to the restriction map $f$, we obtain a finite set $\Sigma \subset C$ such that for every point $y\in C\backslash \Sigma$, all irreducible components of $f^{-1}(y)$ are algebraically equivalent. Let $E\in |A|$ be a general divisor such that the intersection $E\cap C$ avoids $\Sigma$ and $E\cap C$ supports on finite points. This possible since $A=\mathcal{O}(1)_{|Y}$. The intersection $E\cap C$ gives a divisor $P=\sum_{i=1} ^l p_i$ on $C$.
    By Lemma \ref{fiberclass}, all irreducible components of each $f^{-1}(p_i)$ are algebraically equivalent. Since algebraic equivalence implies homological equivalence, each component $W$ of $f^{-1}P\cap D$ represents a multiple of the homology class $[f^{-1}P\cap D]_{D,hom} = [L]_{|D}\cap [D]_{D,hom}$ in $H_{2n-4}(D,\mathbb{R})$.

    This completes the proof.
\end{proof}

\section{Proof of the main result}\label{sec main result}

In this section, we give a self-contained proof of Theorem \ref{mainthrm}.

\subsection{Baby case}
We start with the baby case when $\dim X =3$.

\begin{prop}\label{n=3}
    Let $X$ be a smooth projective variety of dimension $3$ and let $L$ be a big and free class on $X$. Then for the operator $L: H^{1,1}(X)\rightarrow H^{2,2}(X)$, we have
    $\ker L = V_{L, \eff}$.
\end{prop}


Before giving the proof, we first fix an ample class $A$ on $X$. We emphasis that we are working with $H^{1,1}(X)$, thus there is some subtlety to deal with involving  movable classes.

\begin{lem}\label{mov3}
Setting as in Proposition \ref{n=3}, then for any nonzero movable class $[D]$, $L\cdot [D]\neq 0$.
\end{lem}

\begin{proof}
If $[D]$ is a movable class, then for any nef class $F$, $F\cdot [D]$ is in the dual of $\Psef^1 (X)$. This can be checked as follows. By definition, $[D]=\lim_k \pi_{k*} \widehat{\omega}_k$ for a sequence of modifications $\pi_k: \widehat{X}_k\rightarrow X$ and K\"ahler classes $\widehat{\omega}_k$ on $\widehat{X}_k$, thus for any $T\in \Psef^1 (X)$,
\begin{equation*}
  F\cdot [D]\cdot T = \lim_k \pi_k ^* F \cdot \widehat{\omega}_k \cdot \pi_k ^* T \geq 0.
\end{equation*}

When $L\cdot [D]=0$, we have $L\cdot (A\cdot [D])=0$. As $A\cdot [D]$ is in the dual of $\Psef^1 (X)$ and $L$ is an interior point of $\Psef^1 (X)$, we must have $A\cdot [D]=0$, yielding that $[D]=0$.

\end{proof}

\begin{lem}\label{mov31}
Let $X$ be a smooth projective variety of dimension $3$. For a prime divisor class or a nonzero movable class $[D]$ on $X$, the quadratic form on $H^{1,1}(X)$ given by $Q(x,y):=[D]\cdot x\cdot y$ has exactly one positive eigenvalue.
\end{lem}

\begin{proof}
It is clear when $D$ is a prime divisor. 

If $[D]$ is a nonzero movable class, as in Lemma \ref{mov3}, we write $[D]=\lim_k \pi_{k*} \widehat{\omega}_k$. The result follows, since $Q(A, A)>0$, 
\begin{equation*}
  Q(x,y)=[D]\cdot x\cdot y = \lim_k (\widehat{\omega}_k \cdot \pi_k ^*x\cdot \pi_k ^* y),
\end{equation*}
the quadratic form defined by $\widehat{\omega}_k$ has exactly one positive eigenvalue and $\pi_k ^*$ is injective on $H^{1,1}(X)$.
\end{proof}

Now we commence the proof of Proposition \ref{n=3}.

\begin{proof}[Proof of Proposition \ref{n=3}]
By definition $V_{L, \eff} \subset \ker L$, thus we only need to show that any $\alpha\in \ker L$ can be written as a combination of divisors annihilated by $L$.

Fix $\alpha\in \ker L$. In general, one cannot expect any positivity from $\alpha$. We prove that one can obtain positivity from the square $\alpha^2$ by modifying $\alpha$ using $V_{L, \eff}$. 

\begin{claim}
There exists $\beta$ such that $\beta -\alpha \in V_{L, \eff}$ and $-\beta^2$ is in the dual of the pseudo-effective cone $\Psef^1 (X)$.
\end{claim}

We first note that the claim implies $\alpha \in V_{L, \eff}$ as desired. Assume the above claim holds. As $L \cdot \beta =0$, we have $-\beta^2 \cdot L = 0$. Combining with the bigness of $L$ (i.e., $L$ lying in the interior of $\Psef^1 (X)$) and that $-\beta^2$ is in the dual of the pseudo-effective cone $\Psef^1 (X)$, we must have $\beta^2=0$. This yields that 
\begin{equation*}
  0 = (A\cdot L \cdot \beta)^2 =(A\cdot L^2 )(A\cdot \beta^2),
\end{equation*}
which by Lemma \ref{equalAF} shows $A\cdot \beta =0$. Thus, $\beta=0$ by the hard Lefschetz property of $A$, implying $\alpha \in V_{L, \eff}$.

Next we prove the existence of $\beta$ as claimed.

For a prime divisor class or a nonzero movable class $[D]$ on $X$, by $\alpha \in \ker L$ we have $L\cdot \alpha\cdot [D]=0$. Then
\begin{equation*}
  0= (L \cdot \alpha\cdot [D])^2 \geq (L ^2 \cdot [D])(\alpha ^2 \cdot [D]),
\end{equation*}
where the inequality follows from Lemma \ref{mov31}.
Hence,
\begin{equation}\label{mov}
  \alpha^2 \cdot [D]  \leq 0
\end{equation}
whenever $L^2 \cdot [D]>0$. Indeed, (\ref{mov}) still holds whenever the weaker positivity $L\cdot [D]\neq 0$ holds. For those $[D]$ with $\alpha ^2 \cdot [D]  \leq 0$, there is nothing to prove. Thus we only need to consider the case when $\alpha ^2 \cdot [D] \geq 0$. Then by applying Proposition \ref{proportional}, there is some real number $c$ such that 
\begin{equation}\label{prop1}
\alpha\cdot [D] = c L\cdot [D].
\end{equation}
Multiplying $\alpha$ on both sides of (\ref{prop1}) yields $\alpha ^2 \cdot [D] =0$. As any modification $\beta$ of $\alpha$ via $V_{L, \eff}$ is still in $\ker L$, what we have proved is that: for any modification $\beta$ of $\alpha$ via $V_{L, \eff}$,
\begin{equation}\label{mov0}
 -\beta^2 \cdot [D] \geq 0 \text{ for any prime divisor class} [D] \text{ with } L\cdot [D]\neq 0 \text{ and any nonzero movable class } [D]. 
\end{equation}

To check that $-\beta^2$ is in the dual of $\Psef^1(X)$, via Lemma \ref{zariski},
it remains to consider the prime divisor classes such that $L\cdot [D]= 0$.  This is the point where we need the modification via $V_{L, \eff}$.

As $L$ is free, the prime divisors $D$ such that $L\cdot [D]= 0$ are exactly the hypersurfaces contracted to points by the morphism 
\begin{equation*}
 \varphi_L : X \rightarrow Y_L,
\end{equation*}
where $\varphi_L$ the Kodaira map given by $H^0(X, L)$. Since $L$ is big, equivalently, $\dim Y_L =3$, the set of the prime divisors $D$ contracted by $\varphi_L$ is finite, which we denote it by $$\mathcal{D}_1=\{D_1,..., D_m\}.$$

Next we prove that up to some modification of $\alpha$ using an element in $V_{L, \eff}$, we can find $\beta$ such that 
\begin{equation}\label{mov01}
\beta^2 \cdot [D_i] \leq 0 \text{ for any } i\in [m].
\end{equation}
Combining with (\ref{mov0}), $-\beta^2 \cdot [E] \geq 0$ for any pseudo-effective class $E$. This shall finish the proof of the existence of the claimed $\beta$.
To this end, by using Lemma \ref{zariski}, we can add finite prime divisor classes and finite movable classes (as we deal with $H^{1,1}(X)$), which we denote by $$\mathcal{D}_2=\{[D_{m+1}], ..., [D_N]\}$$ 
to the set $\mathcal{D}_1$, such that $\mathcal{C}=\mathcal{D}_1\cup \mathcal{D}_2$ spans $H^{1,1}(X)$. 
As $\mathcal{C}$ spans the whole space, we can write 
\begin{equation}\label{span1}
 \alpha=\sum_{i=1} ^m a_i [D_i] + \sum_{j=m+1}^N a_j [D_j],
\end{equation}
for some real numbers $a_1,...,a_N$.
As $L\cdot \alpha=0$ and $L\cdot [D_i]=0$ for $i\in [m]$, we indeed have $$\sum_{j=m+1}^N a_j [D_j] \in \ker L.$$ 

\begin{claim}
There exist real numbers $b_1,...,b_m$ such that the class 
\begin{equation*}
\beta:= \sum_{i=1} ^m b_i [D_i] + \sum_{j=m+1}^N a_j [D_j]
\end{equation*}
satisfies that for any $k\in [m]$,
\begin{equation}\label{bcoef}
 A\cdot [D_k] \cdot \beta =0.
\end{equation}
\end{claim}

By applying Hodge index theorem on the prime divisor $D_k$, $1\leq k\leq m$, to $$A\cdot [D_k] \cdot \beta =0$$ as claimed, we get $\beta^2 \cdot [D] \leq 0$ for any prime divisor $D$ with $L\cdot [D]=0$. This is exactly what we want in (\ref{mov01}).

Now we prove the existence of $b_1,...,b_m$ such that (\ref{bcoef}) holds. We endow $\mathbb{R}^{[N]}$ with the Euclidean inner product $\langle-, -\rangle$. Consider the linear operator $\mathcal{A}$ defined by the $N \times N$ symmetric matrix $$\mathcal{A} = [A\cdot D_i \cdot D_j]_{1\leq i, j\leq N}$$ and the projection maps $$p_1 (x_1,..,x_m, x_{m+1},...,x_N)^T = (0,..,0, x_{m+1},...,x_N)^T$$ and $$p_2 (x_1,..,x_m, x_{m+1},...,x_N)^T = (x_1,..,x_m,0,..,0)^T.$$ Then the existence of $b_1,...,b_m$ is equivalent to the solvability of the linear system
\begin{equation}\label{eqlin}
 p_2 \circ \mathcal{A} \left[
                         \begin{array}{c}
                           x_1 \\
                           \vdots \\
                           x_m \\
                           a_{m+1} \\
                           \vdots \\
                           a_N \\
                         \end{array}
                       \right] =0.
\end{equation}
Denote $\mathbf{x} = (x_1,...,x_m, a_{m+1},...,a_N)^T$, then the equation (\ref{eqlin}) can be rewritten as 
\begin{equation}\label{eqlin1}
 p_2 \circ \mathcal{A}\circ p_2 (\mathbf{x})
                       =
 -p_2 \circ \mathcal{A}\circ p_1 (\mathbf{x}) =-p_2 \circ \mathcal{A} \left[
                         \begin{array}{c}
                           0 \\
                           \vdots \\
                           0 \\
                           a_{m+1} \\
                           \vdots \\
                           a_N \\
                         \end{array}
                       \right]. 
\end{equation}
By (\ref{eqlin1}), the solvability is equivalent to $$\Image (p_2 \circ \mathcal{A}\circ p_2) \supset \Image (p_2 \circ \mathcal{A}\circ p_1),$$
which, by taking orthogonal complements in $\mathbb{R}^{[N]}$ with respect to the Euclidean inner product $\langle-, -\rangle$ and using the self-adjoint property of $p_2 \circ \mathcal{A}\circ p_2$, is equivalent to 
\begin{equation}\label{kercont}
\ker (p_2 \circ \mathcal{A}\circ p_2) \subset \Image (p_2 \circ \mathcal{A}\circ p_1)^{\perp}.
\end{equation}

We prove (\ref{kercont}) as follows. If the vector $\mathbf{y}=(y_1,...,y_m,y_{m+1},...,y_{N})^T$ is in $\ker (p_2 \circ \mathcal{A}\circ p_2)$, then we must have $$\sum_{1\leq i,j\leq m}\mathcal{A}_{ij} y_i y_j=0,$$ which by the definition of $\mathcal{A}$ is equivalent to the following intersection number equality:
\begin{equation}\label{h1}
A \cdot \gamma^2 =0,
\end{equation}
where $\gamma:= \sum_{i=1}^m y_i [D_i]$.
Note that $L\cdot \gamma =0$, hence
\begin{equation}\label{h2}
A \cdot L\cdot \gamma =0.
\end{equation}
As $L$ is big, we also have
\begin{equation}\label{h3}
A \cdot L^2 >0.
\end{equation}
Using Lemma \ref{equalAF} to (\ref{h1}), (\ref{h2}) and (\ref{h3}) shows $A\cdot \gamma=0$. 
Now for any $\mathbf{z}\in \mathbb{R}^{[N]}$, 
\begin{equation*}
  \langle \mathbf{y}, p_2 \circ \mathcal{A}\circ p_1 (\mathbf{z})\rangle = \sum_{j={m+1}}^N \sum_{i=1}^m (A\cdot D_i \cdot D_j) y_i z_j = A\cdot \gamma \cdot (\sum_{j={m+1}}^N z_j D_j) =0,
\end{equation*}
thus $$\mathbf{y}\perp \Image (p_2 \circ \mathcal{A}\circ p_1).$$ 
This finishes the proof of (\ref{kercont}), hence the proof of (\ref{bcoef}).

This finishes the proof of the lemma.
\end{proof}

The first claim in the above proof is a special case of the local Hodge index inequality proved in \cite[Theorem B]{huxiaoLefcharArxiv}. We expect that Lemma \ref{n=3} admits a more elementary proof via the geometry of morphisms.

\subsection{Beyond the baby case}
Using the baby case on 3-folds, we prove the following generalization to higher dimension. This will be needed in the general case, and it is also an instructive prototype showing how the induction on dimension works.

\begin{prop}\label{weakLef}
Let $X$ be a smooth projective variety of dimension $n\geq 3$. Let $L$ be a free class with $\nd(L)\geq 3$. Then 
    \begin{equation*}
        \ker (L:H^{1,1}(X)\rightarrow H^{2,2}(X))=\Span_{\mathbb{R}}\{[D]:\ D \in \prim(X), L\cdot [D]=0\}.
    \end{equation*}
\end{prop}


\begin{proof}
Note that we only need to prove that any element in the kernel space can be written as the desired linear span.

The case $n=3$ follows from Proposition \ref{n=3}.
Next we proceed by induction on the dimension of the variety. 
Assume that $n\geq 4$ and the result holds for varieties with dimension at most $n-1$. As $L$ is free and $\nd(L)\geq 3$, a general fiber $F$ of $\varphi_L : X \rightarrow Y_L$ has $\dim F \leq n-3$. Hence, the set
\begin{equation*}
  \Null^2 (L)= \overline{\bigcup_{L\cdot [V]=0, \codim V= 2} V}^{\zar}
  \end{equation*}
where $V$ ranges over irreducible subvarieties of codimension 2, is a proper Zariski closed set. Fix a very ample divisor $A$. By Bertini's theorem, we can take a general irreducible smooth hypersurface $H\in |A|$ such that $H$ intersects each component of $\Null^2(L)$ transversally at a general point and the section $H\cap W$ is irreducible.
The second condition can be satisfied since each component of $\Null^2(L)$ has dimension at least 2.

Take $\alpha\in \ker L$, i.e., $L\cdot \alpha=0$. Restricting it to $H$, we get $L_{|H}\cdot \alpha_{|H}=0$. Note that we still have $\nd(L_{|H})\geq 3$, since $H$ is ample. By induction, there exist finite prime divisors $E_i$ of $H$ with 
\begin{equation}\label{ind0}
L_{|H}\cdot [E_i]_H=0
\end{equation}
and real numbers $c_i$ such that
\begin{equation}\label{ind}
  \alpha_{|H} = \sum c_i [E_i]_H.
\end{equation}
Applying $j_{H*}$ to (\ref{ind0}) shows $L\cdot [E_i]_X =0$, yielding that $E_i$ is contained in $\Null^2(L)$. As $E_i \subset H$, by our choice of $H$, $E_i$ cannot be a component of $\Null^2(L)$, hence it must be contained in a component $D_i$ of $\Null(L)$ with $\dim D_i =n-1$. As $H\cap D_i$ is irreducible and $\dim H\cap D_i =n-2$, we get $E_i = H\cap D_i$. This implies that (\ref{ind}) can be written as
\begin{equation}\label{ind1}
  \alpha_{|H} = \sum c_i [D_i]_{|H}.
\end{equation}
Applying Lefschetz hyperplane theorem to (\ref{ind1}) tells
\begin{equation}\label{ind2}
  \alpha = \sum c_i [D_i].
\end{equation}

We claim that $L \cdot [D_i]=0$, which will finish the proof. To this end, we just need to note that $$0=L_{|H}\cdot [E_i]_H= L_{|H}\cdot [D_i]_{|H}.$$ Thus, $A^{n-2}\cdot L\cdot  [D_i]=0$, implying $L\cdot [D_i]=0$.

This finishes the proof.
\end{proof}

\begin{rmk}
Propositions \ref{n=3} and \ref{weakLef} can be also derived straightforward from \cite[Theorem A]{huxiaoLefcharArxiv} which holds more generally for nef collections.
We provide the proofs here due to two reasons: firstly to make the proof self-contained, and secondly to clarify the simplifications due to freeness in some geometric arguments.  
\end{rmk}

\subsection{The general case}

Recall that we aim to prove the following result:

\begin{thrm}
Let $X$ be a smooth projective variety of dimension $n$. Let $\mathcal{L}=(L_1,...,L_{n-2})$ be a supercritical collection of free classes on $X$.
Then $$\ker \mathbb{L} = V_{\mathcal{L},\eff},$$
where $V_{\mathcal{L},\eff}$ is the vector space generated by prime divisors annihilated by $\mathbb{L}$:
\begin{align*}
V_{\mathcal{L}, \eff} = \Span_{\mathbb{R}} \{[D]: D\in \prim(X),\ \mathbb{L} \cdot [D]=0\}.
\end{align*}
\end{thrm}

\begin{proof}
It is obvious that $V_{\mathcal{L}, \eff}\subset \ker \mathbb{L}$. 
For the converse inclusion, we proceed by induction on the dimension $n$.

The base case $n=3$ is established in Proposition \ref{n=3}. 

In the sequel, we assume that $n\geq 4$ and the result holds for all smooth projective varieties of dimension at most $n-1$. Recall from Lemma \ref{locus} that $\Null^2 (\mathcal{L})$ is a proper Zariski closed subset.
As $L_{n-2}$ is free,
by Bertini's theorem and Lemma \ref{algequiv}, we can choose a general hypersurface with $H \in |L_{n-2}|$ satisfying the following properties:
    \begin{enumerate}
        \item[(a)] $H$ is a smooth irreducible hypersurface of $X$;
        \item[(b)] $H$ meets each irreducible component $D$ of $\Null^2(\mathcal{L})$ transversally at a general point;
        \item[(c)] Let $W$ be an irreducible component of $\Null^2(\mathcal{L})$ and let $Z$ be an irreducible component of $H\cap W$. Then there exists $c\in \mathbb{R}_{>0}$ such that 
        $$[Z]_{W,hom}=c[H\cap W]_{W,hom}=cL_{n-2|W}\cap [W]_{W,hom} \text{ in } H_{2n-4}(W,\mathbb{R}).$$
    \end{enumerate}
    Now let $\alpha\in \ker \mathbb{L}$. Multiplying the equation $\mathbb{L\cdot \alpha}=0$ by $\alpha$ and then by $L_{n-2}$, we obtain
    \begin{equation}\label{eq1}
        \alpha \cdot \alpha \cdot L_1\cdot...\cdot L_{n-3}\cdot L_{n-2}=\alpha_{|H}\cdot \alpha_{|H} \cdot L_{1|H}\cdot...\cdot L_{n-3|H}=0,
    \end{equation}
     and
    \begin{equation}\label{eq2}
        \alpha_{|H}\cdot L_{n-2|H} \cdot L_{1|H}\cdot...\cdot L_{n-3|H}=0.
    \end{equation}
    On the other hand, since $\mathcal{L}$ is supercritical, by Lemma \ref{hallrado} we have
    \begin{equation}\label{eq3}
        L_{n-2|H}\cdot L_{n-2|H} \cdot L_{1|H}\cdot...\cdot L_{n-3|H}>0.
    \end{equation}
    By applying Lemma \ref{equalAF} to (\ref{eq1}), (\ref{eq2}) and (\ref{eq3}), it follows that
    \begin{equation*}
        \alpha_{|H}\cdot L_{1|H}\cdot...\cdot L_{n-3|H}=0.
    \end{equation*}

As $\mathcal{L}$ is supercritical, the restriction collection $(L_{1|H},...,L_{n-3|H})$ is supercritical on the lower dimensional variety $H$. Then by induction there exist irreducible varieties (i.e., prime divisors of $H$) 
    \begin{equation*}
        E_1,...,E_s\subset H, \text{ with } \dim E_i=n-2
    \end{equation*}
     and real coefficients $a_1,...,a_s$ such that in the cohomology group of $H$, we have:
    \begin{equation}\label{2}
        \alpha_{|H} =\sum_{i=1}^s a_i[E_i]_{H}.
    \end{equation}
    Moreover, each $E_i$ satisfies 
    \begin{equation}\label{1}
        L_{1|H}\cdot...\cdot L_{n-3|H}\cdot[E_i]_{H}=0.
    \end{equation}

    Pushing forward (\ref{1}) by $j_{H*}$ shows
    \begin{equation*}
        L_1\cdot...\cdot L_{n-2}\cdot [E_i]_X=0.
    \end{equation*}
    So each $E_i$ lies in the null locus $\Null^2(\mathcal{L})$. Since $E_i\subset H$ and by our choice on $H$, $H$ meets every component of $\Null^2(\mathcal{L})$ properly, thus no $E_i$ can be an irreducible component of $\Null^2(\mathcal{L})$. Hence each $E_i$ is contained in some codimension one component $W_i$ of $\Null^2(\mathcal{L})$. In fact $E_i$ is an irreducible component of $H\cap W_i$, as $E_i$ is itself an irreducible subvariety of dimension $n-2$ and $H$ intersects $W_i$ transversally at a general point. By our choice of $H$, there exists for each $E_i$ a positive constant $c_i\in \mathbb{R}_{>0}$ such that in homology
    \begin{equation}\label{3}
        [E_i]_{W_i,hom}=c_{i} L_{n-2|W_i}\cap [W_i]_{W_i,hom} 
    \end{equation}
    and consequently, by pushing forward via $j_{W_i*}$,
    \begin{equation*}
        [E_i]_{X,hom}=c_{i}L_{n-2}\cap [W_i]_{X,hom}.
    \end{equation*}
    By Poincar\'e duality, this is equivalent to 
    \begin{equation}\label{4}
        [E_i]_{X}=c_{i}L_{n-2}\cdot [W_i]_{X}.
    \end{equation}
    Next, apply $j_{H*}$ to (\ref{2}) and use (\ref{4}), we get
    \begin{equation*}
        L_{n-2}\cdot \alpha =\sum_{i=1}^s a_i [E_i]_X=\sum_{i=1}^s a_i c_i L_{n-2}\cdot [W_i]_{X}=L_{n-2}\cdot \sum_{i=1}^s a_i c_i [W_i]_{X}.
    \end{equation*}
    Hence,
    \begin{equation*}
        L_{n-2}\cdot(\alpha-\sum_{i=1}^s a_i c_i [W_i]_{X})=0.
    \end{equation*}

   By Proposition \ref{weakLef} we conclude
    \begin{equation*}
        \alpha-\sum_{i=1}^s a_i c_i [W_i]_X= \sum_{L_{n-2}\cdot [D]=0}a_D[D]
    \end{equation*}
    where the sum runs over all prime divisors $D\subset X$ with $L_{n-2}\cdot [D]=0$ (in particular, $\mathbb{L}\cdot [D]=0$). 
    The theorem then follows once we verify
    \begin{equation*}
        \mathbb{L}\cdot [W_i]_X=0, \forall i=1,...,s,
    \end{equation*}
    which we check below.
    
    By Poincar\'e duality, $\mathbb{L}\cdot [W_i]_X=0$ is equivalent to
    \begin{equation*}
       ( L_1\cdot...\cdot L_{n-2})\cap [W_i]_{X,hom}=0.
    \end{equation*}
    The projection formula gives 
    \begin{equation*}
        ( L_1\cdot...\cdot L_{n-2})\cap [W_i]_{X,hom}=j_{W_i*}((L_{1|W_i}\cdot...\cdot L_{n-2|W_i})\cap [W_i]_{W_i,hom}).
    \end{equation*}
    By the associativity of cap products, together with (\ref{3}), we obtain
    \begin{align*}
        (L_{1|W_i}\cdot...\cdot L_{n-2|W_i})\cap [W_i]_{W_i,hom}&=(L_{1|W_i}\cdot...\cdot L_{n-3|W_i})\cap (L_{n-2|W_i}\cap [W_i]_{W_i,hom})\\
        &=\frac{1}{c_i} (L_{1|W_i}\cdot...\cdot L_{n-3|W_i})\cap [E_i]_{W_i,hom}.
    \end{align*}
    Again by the projection formula, we have
     \begin{equation*}
         (L_{1|W_i}\cdot...\cdot L_{n-3|W_i})\cap [E_i]_{W_i,hom}=(j_{E_i\subset W_i})_*((L_{1|E_i}\cdot...\cdot L_{n-3|E_i})\cap[E_i]_{E_i,hom}).
     \end{equation*}
     Combining all above equations, we obtain
     \begin{equation*}
         ( L_1\cdot...\cdot L_{n-2})\cap [W_i]_{X,hom}=\frac{1}{c_i} j_{W_i*}(j_{E_i\subset W_i})_*((L_{1|E_i}\cdot...\cdot L_{n-3|E_i})\cap[E_i]_{E_i,hom}).
     \end{equation*}
     Chaining these pushforwards together and using functoriality
     \begin{equation*}
         j_{W_i*}(j_{E_i\subset W_i})_*=(j_{E_i\subset X})_*=j_{H*} (j_{E_i\subset H})_*,
     \end{equation*}
     we obtain
     \begin{align*}
         ( L_1\cdot...\cdot L_{n-2})\cap [W_i]_{X,hom}&=\frac{1}{c_i}j_{H*} (j_{E_i\subset H})_*
         ((L_{1|E_i}\cdot...\cdot L_{n-3|E_i})\cap[E_i]_{E_i,hom})\\
         &=\frac{1}{c_i}j_{H*}((L_{1|H}\cdot...\cdot L_{n-3|H})\cap [E_i]_{H,hom})\\
         &=0,
     \end{align*}
    where the last equality follows from (\ref{1}).

    This finishes the proof.
\end{proof}

\section{Applications}\label{sec appl}

\subsection{The extremals for the AF inequality}
Applying Theorem \ref{mainthrm} and the construction via toric geometry (see e.g. \cite{fultonToricBook}, \cite{coxToricBOOK}), we shall give an algebro-geometric proof of the characterization of the extremals of the AF inequality for a supercritical collection of rational convex polytopes, originally due to \cite[Theorem 8.1 or Corollary 2.16]{handelAFextremACTA}. Theorem \ref{mainthrm} is applicable for rational polytopes, since on a smooth projective toric variety every nef line bundle is semiample.

Given a collection $\mathcal{P}=(P_1,...,P_{n-2}) \subset \mathbb{R}^n$ of convex bodies, we call it supercritical if 
$$\dim (P_I)\geq |I|+2,\ \forall I\subset [n-2].$$
For a convex body $K\subset \mathbb{R}^n$, it is uniquely determined by its support function $h_K$ defined by 
\begin{equation*}
  h_K (u)= \sup_{x\in K} x\cdot u.
\end{equation*}
Given a nonzero vector $u\in \mathbb{R}^n$, we denote by
\begin{equation*}
  F(K, u)=\{x\in K: h_K (u) = u\cdot x\},
\end{equation*}
the unique face of $K$ with outer normal direction $u$. 

The following is clear by the definition.

\begin{lem}\label{suppvec}
For any convex bodies $P, Q$ and $v\in \mathbb{R}^n$,
    \begin{enumerate}
        \item $h_{P+Q}(u)=h_P (u)+h_Q(u)$;
        \item $h_{P+v} (u)=h_P (u)+v\cdot u$;
        \item $F(P+Q, u)=F(P,u)+F(Q,u)$.
    \end{enumerate}
\end{lem}

\begin{defn}
    Let $\mathcal{P}=(P_1,...,P_{n-2})$ be a collection of convex bodies in $\mathbb{R}^n$, and let $B$ be the unit ball in $\mathbb{R}^n$. A unit vector $u\in S^{n-1}$ is called a $(B,\mathcal{P})$-extreme normal direction if $u\in \supp S_{B,\mathcal{P}}$, where $S_{B,\mathcal{P}}$ is the surface area measure given by $B, P_1,...,P_{n-2}$.
    
\end{defn}

We will also refer to any nonzero vector $u\in \mathbb{R}^n\backslash \{0\}$ a $(B,\mathcal{P})$-extreme normal direction if its  normalization $\frac{u}{|u|}$ lies in $\supp S_{B,\mathcal{P}}$. This convention is convenient in toric geometry, where it is common to work with primitive lattice vectors rather than unit vectors.

If $P_1,...,P_{n-2}$ are polytopes, the support of mixed surface area measure $\supp S_{B,\mathcal{P}}$ admits an explicit geometric description due to Shenfeld and van Handel. Let $Q$ be a polytope that contains $P_1,...,P_{n-2}$ as Minkowski summands. Denote by $\{F_1,...,F_m\}$ the set of all facets of $Q$, and let $\{u_1,...,u_m\}$ the their corresponding outer unit normal vectors. Define the set
$$E_P=\{(i,j)\in [m]\times [m]:\dim  F_i\cap F_j=n-2\}$$
consisting all pair of facets that intersect in a codimenion-two face.
For each such pair $(i,j)\in E_P$, the vectors $u_i$ and $u_j$ form an angle strictly between 0 and $\pi$, there is a unique shortest geodesic segment $e_{ij} \subset S^{n-1}$ connecting them. 
By construction, for any $u\in \relint e_{ij}$, the face $F(Q,u)$ has dimension $n-2$. We call the edge $e_{ij}$ active if 
$$V_{n-2}(F(P_1,u),...,F(P_{n-2},u))>0,\forall u \in \relint e_{ij}$$
and denote by $E\subset E_P$ the set of all active edges.
The following characterization of the support of $S_{B,\mathcal{P}}$ is established in \cite[Lemma 5.5]{handelAFextremACTA}.

\begin{lem}\label{extremehallrado}
    Let $\mathcal{P}=(P_1,...,P_{n-2})$ be a collection of polytopes in $\mathbb{R}^n$ and let $u\in S^{n-1}$. The following statements are equivalent:
    \begin{enumerate}
        \item $u\in \supp S_{B,\mathcal{P}}$;
        \item $\dim F(P_I,u)\geq |I|, \forall I\subset [n-2]$;
        \item $u\in e_{ij}$ for some $(i,j)\in E$.
    \end{enumerate}
\end{lem}

\begin{cor}\label{extreme}
    Let $M,N$ be rational polytopes in $\mathbb{R}^n$ and let $Q$ be a rational polytope admitting both $M$ and $N$ as Minkowski summands. Denote by $\{u_1,...,u_m\}$ the set of outer primitive normal vectors of the facets of $Q$.
    Then $M$ and $N$ have the same
    supporting hyperplanes in all $(B, \mathcal{P})$-extreme normal directions if
    \begin{equation*}
        h_M(u_i)=h_N(u_i), \forall \frac{u_i}{|u_i|} \in  \supp S_{B,\mathcal{P}}.
    \end{equation*}
\end{cor}
\begin{proof}
    We aim to show
    \begin{equation*}
        h_M(u)=h_N(u) , \forall u\in \supp S_{B,\mathcal{P}}.
    \end{equation*}
    By Lemma \ref{extremehallrado}, the support decomposes as a union of edges
    \begin{equation*}
        \supp S_{B,\mathcal{P}}=\bigcup_{(i,j)\in E} e_{ij}.
    \end{equation*}
    Fix an arbitrary edge $e_{ij}$ with endpoints $\frac{u_i}{|u_i|},\frac{u_j}{|u_j|}$. By assumption, $h_M$ and $ h_N$ agree at $u_i$ and $u_j$. Since $M,N$ are both Minkowski summands of $Q$, their support functions are linear on the cone generated by 
    $u_i, u_j$ and thus they agree there. Note that the edge $e_{ij}$ lies within this cone as its unit section. It follows that $h_M=h_N$ on the entire edge $e_{ij}$. Since $e_{ij}$ is arbitrary, the equality holds on the whole support, completing the proof.
\end{proof}

\begin{thrm}\label{corhandel1}
Let $\mathcal{P}=(P_1,...,P_{n-2}) \subset \mathbb{R}^n$ be a supercritical collection of rational polytopes. Assume that $M,N$ are rational polytopes such that $V(M, N,\mathcal{P} )>0$. Then 
\begin{equation}\label{afeq}
  V(M, N, \mathcal{P})^2 = V(M, M, \mathcal{P})V(N, N, \mathcal{P})
\end{equation}
if and only if there exist $a\in \mathbb{R}_{>0}$ and $v\in \mathbb{R}^n$ such that $M$ and $aN+v$ have the same
supporting hyperplanes in all $(B, \mathcal{P})$-extreme normal directions. 
\end{thrm}

Before giving the proof, let us recall some standard facts from toric geometry.
Let $Q\subset \mathbb{R}^n$ be a full-dimensional rational polytope. Denote by $X_Q$ the associated projective toric variety. For simplicity, we assume that $X_Q$ is smooth, in which case the polytope is called Delzant. The torus-invariant prime divisors $D_1,...,D_m$ corresponds bijectively to the facets $F_1,...,F_m$ of $Q$, with outer primitive normal vectors denoted by $u_1,...,u_m$. 
If $P$ is a Minkowski summand of $Q$, we define the semiample $\mathbb{Q}$-divisor associated to $P$ as
$$L_P=\sum_{i}h_P(u_i)D_i,$$
where $h_P$ is the support function.
The Bernstein-Khovanskii-Kushnirenko theorem states that for any tuple of rational polytopes $(P_1,...,P_n)$, each a Minkowski summand of $Q$, we have
\begin{equation*}
  V(P_1,...,P_n) = \frac{1}{n!} (L_{P_1} \cdot...\cdot L_{P_n}).
\end{equation*}
As a consequence,
\begin{equation*}
    \nd(L_P)=\dim P.
\end{equation*}

The following result should be well known to experts, we include a proof as we do not find an exact reference. 

\begin{lem}\label{restriction}
    Let $X_Q$ be the smooth projective toric variety associated to a Delzant polytope $Q$, and let $P$ be a rational polytope that is a Minkowski summand of $Q$. Then for any facet $F$ of $Q$ with outer primitive normal vector $u_F$, we have
    \begin{equation*}
        L_{P|D_F} \sim_{\mathbb{Q}} L_{F(P,u_F)}.
    \end{equation*}   
    where $D_F$ is the torus-invariant prime divisor corresponding to $F$, which is itself a smooth projective toric variety associated to the polytope $F$, and $L_{F(P,u_F)}$ is the torus-invariant divisor determined by $F(P,u_F)$ which can be viewed as a Minkowski summand of $F$ up to translation.
\end{lem}
\begin{proof}
    By choosing a suitable basis of the lattice and translating the polytopes if necessary, we may assume that $F \subset \{x_1 = 0\}$ and $Q \subset \{x_1 \leq 0\}$. Thus, the outer primitive normal vector of $Q$ at $F$ becomes $u_F = e_1$, the first standard basis vector. Furthermore, we may assume (up to translation) that the face $F(P, e_1)$ is contained in $F$.
    Let $G_1, \ldots, G_l$ denote the facets of $Q$ that intersect $F$ in codimension two, and let $u_1, \ldots, u_l$ be their respective outer primitive normal vectors. Then the restriction of $L_P$ to $D_F$ is given by
    \begin{equation*}
        L_{P|D_F}=\sum_{i=1}^l h_{P}(u_i) D_i\cap D_F,
    \end{equation*}
    where $D_i$ is the torus-invariant prime divisor corresponding to $G_i$. Here we use the smoothness of $X_Q$ to ensure that $D_i$ and $D_F$ intersect transversely, so that $D_{i|D_F} = D_i \cap D_F$.
    On the other hand, viewing $F(P,u_F)$ as a Minkowski summand of $F$, the corresponding divisor $L_{F(P,u_F)}$ on $D_F$ is given by
    $$L_{F(P,u_F)}=\sum_{i=1}^lh_{F(P,u_F)}(v_i)D_i\cap D_F,$$
    where $v_i$ is the outer primitive normal vector to the facet $G_i \cap F$ of $F$.
    It remains to show that $h_P(u_i) = h_{F(P,u_F)}(v_i)$ for all $i=1,\ldots,l$. To this end, observe that the projection of $u_i$ onto the hyperplane $\{x_1 = 0\}$ is
    \begin{equation*}
        u_i- (u_i\cdot e_1)e_1.
    \end{equation*}
    Since this vector lies in the direction of $v_i$, which is primitive, there exists $\lambda \in \mathbb{Z}_{>0}$ such that
    \begin{equation}\label{461}
        u_i- (u_i\cdot e_1)e_1=\lambda v_i.
    \end{equation}
    We now argue that $\lambda = 1$. Because $Q$ is Delzant, the set $\{u_i, e_1\}$ can be completed to a $\mathbb{Z}$-basis of the lattice. Note (\ref{461}) expresses $v_i$ as
    \begin{equation*}
        v_i= \frac{1}{\lambda}u_i- \frac{u_i\cdot e_1}{\lambda}e_1.
    \end{equation*}
    which is a linear combination of $u_i$ and $e_1$ with rational coefficients. Since $v_i$ is a primitive lattice vector, these coefficients must be integers, forcing $\lambda = 1$.
    
    Finally, let $x\in G_i\cap F \subset \{x_1=0\}$. Then
    \begin{equation*}
        h_P(u_i)=x\cdot u_i=x\cdot (v_i+(u_i\cdot e_1)e_1)=x\cdot v_i=h_{F(P,u_F)}(v_i)
    \end{equation*}
    since $x\cdot e_1=0$. This finishes the proof.
\end{proof}

\begin{lem}\label{hallradotoric}
    Let $X_Q$ be the smooth projective toric variety associated to a Delzant polytope $Q$, and let $P_1,...,P_{n-2}$ be Minkowski summands of $Q$, with the corresponding divisors given by $L_1,...,L_{n-2}$. Then, for any torus-invariant divisor $D_F$ corresponding to a facet $F$ of $Q$, with outer primitive normal vector $u_F$, we have
    \begin{equation*}
        \mathbb{L}\cdot  [D_F]=0 \iff \frac{u_F}{|u_F|}\notin \supp S_{B,\mathcal{P}}.
    \end{equation*}
\end{lem}
\begin{proof}    
    First, we show
    \begin{equation}\label{421}
        L_1\cdot...\cdot L_{n-2}\cdot [D_F]=0 \iff L_{1|D_F}\cdot...\cdot L_{n-2|D_F}=0.
    \end{equation}
    In (\ref{421}), only the $\implies$ direction is nontrivial.
    We fix an ample divisor $A$ on X, the assumption tells
    \begin{equation}\label{423}
        A\cdot L_1\cdot...\cdot L_{n-2}\cdot [D_F]=A_{|D_F}\cdot L_{1|D_F}\cdot...\cdot L_{n-2|D_F}=0.
    \end{equation}
    Since $A_{|D_F}$ is still ample and the $L_{i|D_F}$ are nef, 
    (\ref{423}) forces
    $$L_{1|D_F}\cdot...\cdot L_{n-2|D_F}=0.$$
    This proves (\ref{421}).
    
    The hypersurface $D_F$ itself, as a smooth projective variety, is the toric variety associated to the facet $F$. Moreover, for any $I\subset [n-2]$, $L_{I|D_F}$ is the divisor given by the polytope $F(P_I,u_F)$ which can be viewed as a Minkowski summand of $F$ up to a translation by Lemma \ref{restriction}. Hence, $$\nd(L_{I|D_F})=\dim F(P_I,u_F).$$
    Then by Lemma \ref{hallrado},
    \begin{equation}\label{462}
        L_{1|D_F}\cdot...\cdot L_{n-2|D_F}=0\iff \exists I\subset [n-2] \text{ s.t. }\dim F(P_I,u_F)<|I| .
    \end{equation}
    The result then follows immediately from (\ref{421}), (\ref{462}) and Lemma \ref{extremehallrado}.
\end{proof}

Now we give the proof of Theorem \ref{corhandel1}.

\begin{proof}[Proof of Theorem \ref{corhandel1}]
By \cite[Section 5.4]{fultonToricBook}, we can construct a smooth projective toric variety $X_Q$ associated to the
Delzant polytope $Q$, which admits $M,N,P_1,...,P_{n-2}$ as Minkowski summands. Denote the corresponding divisors by $L_M, L_N, L_{P_1}, ...,L_{P_{n-2}}$, all of which are semiample.

The equality $$ V(M, N, \mathcal{P})^2 = V(M, M, \mathcal{P})V(N, N, \mathcal{P})$$ is equivalent to
    \begin{equation}\label{afeq1}
        (L_M\cdot L_N\cdot L_{P_1}\cdot...\cdot L_{P_{n-2}})^2 = (L_M^2\cdot L_{P_1}\cdot...\cdot L_{P_{n-2}})(L_N^2\cdot L_{P_1}\cdot...\cdot L_{P_{n-2}})
    \end{equation}
    which, according to Lemma \ref{equalAF}, implies that there exists $a\in \mathbb{R}$ such that
    \begin{equation*}
        (L_M-aL_N)\cdot L_{P_1}\cdot...\cdot L_{P_{n-2}}=0.
    \end{equation*}
    Applying Theorem \ref{mainthrm}, we obtain 
    \begin{equation}\label{411}
        L_M-aL_N= \sum a_D[D] \text{ in } H^{1,1}(X)
    \end{equation}
    where the sum runs over all prime divisors $D$ with $\mathbb{L}\cdot [D]=0$. Since the collection $(L_{P_1},...,L_{_{P_{n-2}}})$ is supercritical,  Remark \ref{rigid} shows that each $[D]$ spans an extremal ray of the pseudo-effective cone and hence must be torus-invariant. Thus
    \begin{equation*}
        L_M-aL_N= \sum_{\mathbb{L}\cdot [D_i]=0} a_i[D_i] \text{ in } H^{1,1}(X)
    \end{equation*}
    where $D_1,...,D_m$ are all torus-invariant prime divisors. Every $D_i$ corresponds to a facet $F_i$ of $Q$. We denote the outer primitive normal vector of $F_i$ by $u_i$.
    
    On a smooth toric variety, rational equivalence and cohomological equivalence coincide by \cite{fulton_sturmfels}, so 
    \begin{equation*}
        L_M-aL_N-\sum_{i}a_i D_i \sim_{\mathbb{R}}0.
    \end{equation*}
    Moreover, torus-invariant divisors that are rationally equivalent to zero can be written as $\mathbb{R}$-linear combinations of principal divisors given by characters $\{ div (\chi^{m})\}$ according to \cite{fulton_sturmfels}. Hence, there is a vector $v\in \mathbb{R}^n$ such that
    \begin{equation*}
        L_M-aL_N-\sum_{\mathbb{L}\cdot [D_i]=0}a_i D_i=\sum_{i}(u_i\cdot v) D_i.
    \end{equation*}

    By Lemma \ref{hallradotoric}, $\mathbb{L}\cdot [D_i]=0$ is equivalent to $\frac{u_i}{|u_i|} \notin \supp S_{B,\mathcal{P}}$. Comparing the coefficients of $D_i$ for those $\frac{u_i}{|u_i|}\in \supp S_{B,\mathcal{P}}$ gives
    \begin{equation*}
        h_M(u_i)-ah_N(u_i)=u_i\cdot v ,\forall \frac{u_i}{|u_i|}\in \supp S_{B,\mathcal{P}}
    \end{equation*}
    Finally, by Lemma \ref{suppvec} and Corollary \ref{extreme}, these equalities hold iff
    \begin{equation*}
        h_M(u)=h_{aN+v}(u) ,\forall u\in \supp S_{B,\mathcal{P}}.
    \end{equation*}

    This finishes the proof.
\end{proof}

\subsection{The extremals for the Khovanskii-Teissier inequality}
 Another consequence is the characterization on the extremals of the Khovanskii-Teissier inequality
 for the intersection numbers of two free divisor classes.
 
 \begin{cor}\label{corlogc1}
Let $X$ be a smooth projective variety of dimension $n$, and let $A, B$ be free divisor classes on $X$. Fix $1\leq k\leq n-1$. Assume that $(A^k \cdot B^{n-k})>0$,  then 
\begin{equation*}
  (A^k \cdot B^{n-k})^2 =  (A^{k-1} \cdot B^{n-k+1}) (A^{k+1} \cdot B^{n-k-1})
\end{equation*}
if and only if there exist some $c>0$ and $c_i\in \mathbb{R}$ such that
\begin{equation*}
  A-cB= \sum_{i=1}^N c_i [D_i]
\end{equation*}
where the $D_i$ are primes divisors such that $A^{k-1} \cdot B^{n-k-1}\cdot [D_i] =0$.
\end{cor}

\begin{proof}
We may assume $k\leq n-k$.

When $k=1$, by the condition $(A^k \cdot B^{n-k})>0$ and the equality, $B$ must be big.

When $k\geq 2$,  by Lemma \ref{hallrado}, the conditions $$(A^{k+1} \cdot B^{n-k-1})>0, (A^{k-1} \cdot B^{n-k+1})>0$$ hold iff $$\nd(A)\geq k+1 (\geq 3), \nd(B)\geq n-k+1 (\geq 3), \nd(A+B)=n.$$

In any case, the collection $\{A[k-1], B[n-k-1]\}$, with $A$ appearing $k-1$ times and $B$ appearing $n-k-1$ times, is supercritical. Then Theorem \ref{mainthrm} applies.
\end{proof}

As mentioned in Remark \ref{rigid}, the number of such $D_i$ is at most the Picard number of $X$, and these divisors are contained in the augmented base locus of a big class and have very strong rigidity, and in particular, span extremal rays of the pseudo-effective cone $\Psef^1(X)$.

\begin{rmk}
If $(A^k \cdot B^{n-k})=0$, then the equality automatically holds. By Lemma \ref{hallrado}, $(A^k \cdot B^{n-k})=0$ iff either $\nd(A)<k$, or $\nd(B)<n-k$, or $\nd(A+B)<n$.
\end{rmk}

\section*{Acknowledgements}
This work is supported by the National Key Research and Development Program of China (No. 2021YFA1002300) and National Natural Science Foundation of China (No. 11901336). In the August of 2024, the second named author presented our another potential approach to Jie Liu, while it seems that our old approach still has some subtlety, we would like to thank him for helpful discussions.

\bibliography{reference}

\newcommand{\etalchar}[1]{$^{#1}$}
\providecommand{\bysame}{\leavevmode\hbox to3em{\hrulefill}\thinspace}
\providecommand{\MR}{\relax\ifhmode\unskip\space\fi MR }
\providecommand{\MRhref}[2]{%
  \href{http://www.ams.org/mathscinet-getitem?mr=#1}{#2}
}
\providecommand{\href}[2]{#2}
\begin{thebibliography}{HHM{\etalchar{+}}25}

\bibitem[Ale38]{alexandroff1938theorie}
Alexander Alexandrov, \emph{Zur theorie der gemischten volumina von konvexen
  {K}\"orpern. {IV}. {D}ie gemischten {D}iskriminanten und die gemischten
  volumina}, Matematicheskii Sbornik \textbf{45} (1938), no.~2, 227--251.

\bibitem[Bou04]{Bou04}
S{\'e}bastien Boucksom, \emph{Divisorial {Z}ariski decompositions on compact
  complex manifolds}, Ann. Sci. \'Ecole Norm. Sup. (4) \textbf{37} (2004),
  no.~1, 45--76. \MR{2050205 (2005i:32018)}

\bibitem[CLS11]{coxToricBOOK}
David~A. Cox, John~B. Little, and Henry~K. Schenck, \emph{Toric varieties},
  Graduate Studies in Mathematics, vol. 124, American Mathematical Society,
  Providence, RI, 2011. \MR{2810322}

\bibitem[CT15]{tosattiNullLocus}
Tristan~C. Collins and Valentino Tosatti, \emph{K\"{a}hler currents and null
  loci}, Invent. Math. \textbf{202} (2015), no.~3, 1167--1198. \MR{3425388}

\bibitem[dCM02]{decataldoLefsemismall}
Mark Andrea~A. de~Cataldo and Luca Migliorini, \emph{The hard {L}efschetz
  theorem and the topology of semismall maps}, Ann. Sci. \'{E}cole Norm. Sup.
  (4) \textbf{35} (2002), no.~5, 759--772. \MR{1951443}

\bibitem[Dem12]{dem_analyticAG}
Jean-Pierre Demailly, \emph{Analytic methods in algebraic geometry}, Surveys of
  Modern Mathematics, vol.~1, International Press, Somerville, MA; Higher
  Education Press, Beijing, 2012. \MR{2978333}

\bibitem[DJV13]{voisinPushforw}
Olivier Debarre, Zhi Jiang, and Claire Voisin, \emph{Pseudo-effective classes
  and pushforwards}, Pure Appl. Math. Q. \textbf{9} (2013), no.~4, 643--664.
  \MR{3263971}

\bibitem[FS97]{fulton_sturmfels}
William Fulton and Bernd Sturmfels, \emph{Intersection theory on toric
  varieties}, Topology \textbf{36} (1997), no.~2, 335--353. \MR{1415592}

\bibitem[Ful93]{fultonToricBook}
William Fulton, \emph{Introduction to toric varieties}, Annals of Mathematics
  Studies, vol. 131, Princeton University Press, Princeton, NJ, 1993, The
  William H. Roever Lectures in Geometry. \MR{1234037}

\bibitem[HHM{\etalchar{+}}25]{huh2025arxiv-realizationshomologyclassesprojection}
Daoji Huang, June Huh, Mateusz Michalek, Botong Wang, and Shouda Wang,
  \emph{Realizations of homology classes and projection areas}, 2025.

\bibitem[HSX23]{hushangxiao2023hardlefArxiv}
Jiajun Hu, Shijie Shang, and Jian Xiao, \emph{Hard {L}efschetz theorems for
  free line bundles}, arXiv:2305.19085 (2023).

\bibitem[Huh22]{huh2022combinatorics}
June Huh, \emph{Combinatorics and {H}odge theory}, Proceedings of the
  {I}nternational {C}ongress of {M}athematicians, 2022.

\bibitem[HX22]{huxiaohardlef2022Arxiv}
Jiajun Hu and Jian Xiao, \emph{Hard {L}efschetz properties, complete
  intersections and numerical dimensions}, arXiv:2212.13548 (2022).

\bibitem[HX23]{huxiaoLefcharArxiv}
Jiajun Hu and Jian Xiao, \emph{Numerical characterization of the hard
  {L}efschetz classes of dimension two}, arXiv:2309.05008 (2023).

\bibitem[Laz04]{lazarsfeldPosI}
Robert Lazarsfeld, \emph{Positivity in algebraic geometry. {I}}, Ergebnisse der
  Mathematik und ihrer Grenzgebiete. 3. Folge. A Series of Modern Surveys in
  Mathematics [Results in Mathematics and Related Areas. 3rd Series. A Series
  of Modern Surveys in Mathematics], vol.~48, Springer-Verlag, Berlin, 2004,
  Classical setting: line bundles and linear series. \MR{2095471}

\bibitem[Pan85]{Panov1987ONSP}
A.~A. Panov, \emph{Some properties of mixed discriminants}, Mat. Sb. (N.S.)
  \textbf{128(170)} (1985), no.~3, 291--305, 446. \MR{815265}

\bibitem[Sch14]{schneiderBrunnMbook}
Rolf Schneider, \emph{Convex bodies: the {B}runn-{M}inkowski theory}, expanded
  ed., Encyclopedia of Mathematics and its Applications, vol. 151, Cambridge
  University Press, Cambridge, 2014. \MR{3155183}

\bibitem[Sta81]{stanleyAF}
Richard~P. Stanley, \emph{Two combinatorial applications of the
  {A}leksandrov-{F}enchel inequalities}, J. Combin. Theory Ser. A \textbf{31}
  (1981), no.~1, 56--65. \MR{626441}

\bibitem[Sta86]{stanleyposet}
\bysame, \emph{Two poset polytopes}, Discrete Comput. Geom. \textbf{1} (1986),
  no.~1, 9--23. \MR{824105}

\bibitem[Sta89]{stanleyLogConc}
\bysame, \emph{Log-concave and unimodal sequences in algebra, combinatorics,
  and geometry}, Graph theory and its applications: {E}ast and {W}est ({J}inan,
  1986), Ann. New York Acad. Sci., vol. 576, New York Acad. Sci., New York,
  1989, pp.~500--535. \MR{1110850}

\bibitem[SvH23]{handelAFextremACTA}
Yair Shenfeld and Ramon van Handel, \emph{The extremals of the
  {A}lexandrov-{F}enchel inequality for convex polytopes}, Acta Math.
  \textbf{231} (2023), no.~1, 89--204. \MR{4652411}

\end{thebibliography}
\bibliographystyle{amsalpha}

\bigskip

\bigskip

\noindent
\textsc{Tsinghua University, Beijing 100084, China}\\
\noindent
\verb"Email: hujj22@mails.tsinghua.edu.cn"\\
\noindent
\verb"Email: jianxiao@tsinghua.edu.cn"

\end{document}